\documentstyle[11pt,leqno,amsfonts,amssymb]{article}
\newtheorem{definition}{Definition}[section]
\newtheorem{lemma}[definition]{Lemma}
\newtheorem{theorem}[definition]{Theorem}
\newtheorem{proposition}[definition]{Proposition}
\newtheorem{corollary}[definition]{Corollary}
\newtheorem{remark}[definition]{Remark}

\newtheorem{example}[definition]{Example}

 1

\newenvironment{proof}{\noindent{\bf Proof~:}}{\QED\medskip}
\def\QED{\hskip0.1em\hfill\null\ \null\nobreak\hfill
\kern3pt\lower1.8pt\vbox{\hrule\hbox
{\vrule\kern1pt\vbox{\kern1.7pt \hbox{$\scriptstyle
QED$}\kern0.2pt}\kern1pt\vrule}\hrule}}

\newcommand{\G}{\Gamma}

\newcommand{\ZZ}{{\mathbb{Z}}}
\newcommand{\CC}{{\mathbb{C}}}

\newcommand{\RR}{{\mathbb{R}}}

\renewcommand{\Im}{{\mathrm{Im}}\,}

\newcommand{\Id}{\hbox{Id}}
\newcommand{\la}{\langle}
\newcommand{\ra}{\rangle}

\begin{document}

\title{\bf The $d\delta$--lemma for
weakly Lefschetz symplectic~manifolds}

\author{Marisa Fern\'andez, Vicente Mu\~noz and Luis Ugarte}
\date{December 20, 2004}

\maketitle

\begin{abstract}
For a symplectic manifold $(M,\omega)$, not necessarily hard
Lefschetz, we prove a version of the Merkulov $d\delta$--lemma
(\cite{Mer, Cav2}). We also study the $d\delta$--lemma and related
cohomologies for compact symplectic solvmanifolds.
\end{abstract}

\section{Introduction} \label{introduction}

Let $(M,\omega)$ be a symplectic manifold, that is, $M$ is a
differentiable manifold of dimension~$2n$ with a closed
non-degenerate $2$--form $\omega$, the {\it symplectic form}.
Denote by $\Omega^k(M)$ the space of the differential $k$--forms
on $M$. According to Libermann \cite{Lib} and Brylinski \cite{Bry}
there is a {\em symplectic star operator} $*:\Omega^k(M)
\longrightarrow \Omega^{2n-k}(M)$ associated to the symplectic
form $\omega$ satisfying $*^2=Id$ (see Section \ref{someresults}
for the definition). Such an operator is the symplectic analogue
of the Hodge star operator on oriented Riemannian manifolds. Then,
one can define the {\em codifferential} $\delta = \pm *d*$ which
satisfies $\delta^2 = 0$ and $d\delta+\delta d=0$ (although
$\delta$ does not satisfy a Leibniz rule \cite{LS}).

As in Riemannian Hodge theory, a $k$--form $\alpha \in
\Omega^k(M)$ is said to be {\em coclosed} if $\delta\alpha=0$,
{\em coexact} if $ \alpha=\delta\beta$ for some $\beta$, and {\em
symplectically harmonic} if it is closed and coclosed. But, unlike
the case of Riemannian manifolds, there are many symplectically
harmonic forms which are exact. This is the reason for which for
any $k\geq 0$, we define the space of {\em harmonic cohomology}
$H^{k}_{\rm hr}(M,\omega)$ of degree $k$ to be the subspace of the
de Rham cohomology $H^k(M)$ consisting of all cohomology classes
which contain at least one symplectically harmonic $k$--form.

Mathieu \cite{Mat} and, independently,  Yan \cite{Ya} proved that
$H^{k}_{\rm hr}(M,\omega)=H^k(M)$ for all $k$ if and only if
$(M,\omega)$ satisfies the hard Lefschetz property, i.e. the map
 $$
  L^{n-k} \colon H^k(M) \longrightarrow H^{2n-k}(M)
 $$
given by $L^{n-k}[\alpha]=[\omega^{n-k}\wedge\alpha]$ is a
surjection for all $k \leq n-1$. On the other hand, for {\em
compact\/} symplectic manifolds, Merkulov and Cavalcanti
(\cite{Mer, Cav2}) showed that the existence of symplectic
harmonic forms in every de Rham cohomology class is equivalent to
the symplectic $d\delta$--{\em lemma}, that is, to the identities
 \begin{equation} \label{eqn:v1a}
 \Im d\cap \ker \delta  = \Im d \delta =   \Im\delta \cap \ker d,
 \end{equation}
which mean that if $\alpha$ is a symplectically harmonic $k$--form
and either is exact or coexact, then $\alpha=d\delta\beta$ for
some $k$--form $\beta$.

Consider the subcomplex  $(\Omega^*_{\delta}(M, \omega),d)$ of the
de Rham complex $(\Omega^*(M),d)$ of $M$, where
$\Omega^k_{\delta}(M,\omega)$  is the space of the coclosed
$k$--forms. We denote by $H^*_{\delta}(M,\omega)$ its cohomology
and by $i$ the natural map
 \begin{equation} \label{eqn:v1aa}
 i\colon H^k_{\delta}(M,\omega)\longrightarrow H^{k}(M),
 \end{equation}
for all $k\geq 0$. In \cite{Gui} Guillemin proved that if $M$ is
compact, then the map $i$ is bijective if and only if $(M,\omega)$
is hard Lefschetz or, equivalently, it satisfies the
$d\delta$--lemma.

In this paper, we aim to generalize these results to symplectic
manifolds which are not hard Lefschetz. Recall the following
definition from \cite{FM}.

\begin{definition} \label{weLefschetz}
A symplectic manifold $(M,\omega)$ of dimension $2n$ is said to be
{\em $s$--Lefschetz\/}, where $0 \leq s\leq n-1$, if the map
 $$
 L^{n-k} \colon H^k(M) \longrightarrow H^{2n-k}(M)
 $$
is an epimorphism for all $k \leq s$. (If $M$ is compact, then we
actually have that $L^{n-k}$ are isomorphisms because of
Poincar\'e duality.)
\end{definition}

Whenever $(M,\omega)$ is not hard Lefschetz, there is some integer
number $s\geq 0$ such that $(M,\omega)$ is $s$--Lefschetz, but not
$(s+1)$--Lefschetz. Note that $(M,\omega)$ is $(n-1)$--Lefschetz
if it satisfies the hard Lefschetz theorem.

Concerning the harmonic cohomology for such manifolds, we have the
following result.

\begin{theorem}\hskip-.1cm{\rm \cite{FMU}}  \label{grados+2}
 Let $(M,\omega)$ be a symplectic manifold of dimension
 $2n$ and let $s\leq n-1$. Then the following statements are
 equivalent:
 \begin{enumerate}
 \item[{\rm (i)}] $(M,\omega)$ is $s$--Lefschetz.
\item[{\rm (ii)}] $H^{k}_{\rm hr}(M,\omega) = H^k(M)$ for every
 $k\leq s+2$, and $H^{2n-k}_{\rm hr}(M,\omega) = H^{2n-k}(M)$ for every $k\leq s$.
\item[{\rm (iii)}] $H^{2n-k}_{\rm hr}(M,\omega) = H^{2n-k}(M)$ for
every $k\leq s$.
 \end{enumerate}
\end{theorem}

Notice that Theorem \ref{grados+2} implies that every de Rham
cohomology class of $M$ admits a symplectically harmonic
representative if and only if $(M,\omega)$ is hard Lefschetz,
which is the result proved independently by Mathieu and Yan
\cite{Mat, Ya}.

For any non-hard Lefschetz symplectic manifold, it seems
interesting to understand how the level $s$ at which the Lefschetz
property is lost affects to other properties of the manifold, such
as the above mentioned $d\delta$--lemma, or to the properties of
the map $i$. Our purpose in this paper is to explore these
questions, as we explain below.

\medskip

In Section \ref{someresults} we recall some properties of the
spaces of harmonic cohomology. In Section \ref{ddeltalemma}, we
sharpen the result of Merkulov and the result of Guillemin by
using the concept of $s$--Lefschetz property. We need first to
weaken the condition of the $d\delta$--lemma to the following

\begin{definition} \label{def:ddeltas}
Let $(M,\omega)$ be a symplectic manifold of dimension $2n$, and
$0 \leq s\leq n-1$. We say that $(M,\omega)$ satisfies the {\em
$d\delta$--lemma up to degree $s$} if
    \begin{equation}\label{eqn:true1}
 \begin{array}{ll}
 \Im d\cap \ker \delta  = \Im d \delta =   \Im\delta \cap \ker d,
  \qquad   & {\mathrm on }\ \Omega^k(M),\ {\mathrm for }\ k \leq s, \\
 \Im d\cap \ker \delta =\Im d \delta,
  \qquad   & {\mathrm on }\ \Omega^{s+1}(M).
 \end{array}
    \end{equation}
\end{definition}

Therefore, if $(M,\omega)$ satisfies the $d\delta$--lemma up to
degree $s$, and $\alpha \in \Omega^{\leq s}(M)$ is symplectically
harmonic and either is exact or coexact then $\alpha=d\delta
\beta$ for some $\beta$; moreover, if $\alpha \in \Omega^{s+1}(M)$
is symplectically harmonic and exact then $\alpha=d\delta \beta$
for some $\beta$.

Following the approach in Cavalcanti's proof~\cite{Cav2} of the
result of Merkulov we prove the following theorem.

\begin{theorem} \hskip-.1cm ($d\delta$--lemma for weakly
Lefschetz manifolds).
 \label{gradosddelta}
 Let $(M,\omega)$ be a compact symplectic manifold of dimension
 $2n$ and let $s\leq n-1$. Then the following statements are
 equivalent:
 \begin{enumerate}
 \item[{\rm (i)}] $(M,\omega)$ is $s$--Lefschetz.
 \item[{\rm (ii)}]  $(M,\omega)$ satisfies the $d\delta$--lemma
 up to degree $s$.
 \item[{\rm (iii)}] The identities {\rm (\ref{eqn:v1a})} hold on
 $\Omega^{\geq (2n-s)}(M)$, and $\Im\delta\cap \ker d =\Im d \delta$ holds
 on $\Omega^{2n-s-1}(M)$.
 \end{enumerate}
\end{theorem}

In Section \ref{ddeltalemma} we also show the following theorem
regarding the map (\ref{eqn:v1aa}) for weakly symplectic
manifolds.

\begin{theorem} \label{mapi}
 Let $(M,\omega)$ be a compact symplectic manifold of dimension
 $2n$ and let $s\leq n-1$. Then the following statements are
 equivalent:
 \begin{enumerate}
 \item[{\rm (i)}] $(M,\omega)$ is $s$--Lefschetz.
 \item[{\rm (ii)}]  The map $i\colon H^k_{\delta}(M,\omega)
 \longrightarrow H^{k}(M)$ is bijective for all
 $k\leq s+1$ and for $k\geq 2n-s$.
 \item[{\rm (iii)}] The map $i\colon H^k_{\delta}(M,\omega)
 \longrightarrow H^{k}(M)$ is bijective for all
 $k\geq 2n-s$.
 \end{enumerate}
\end{theorem}

The harmonic cohomology of compact symplectic nilmanifolds has
been studied by different authors (see \cite{Ya, IRTU, Yam}). In
Section~\ref{solvmanifolds} we consider compact solvmanifolds $M
=\Gamma\backslash G$, where $G$ is a simply connected solvable Lie
group whose Lie algebra ${\frak g}$ is completely solvable, i.e.,
the map ${\rm ad}_X\colon {\frak g} \longrightarrow {\frak g}$ has
only real eigenvalues for any $X \in {\frak g}$, and $\Gamma$ is a
discrete subgroup of $G$ such that the quotient $M
=\Gamma\backslash G$ is compact. We show that the harmonic
cohomology of $(M =\Gamma\backslash G, \omega)$ is isomorphic to
the harmonic cohomology at the level of the invariant forms. We
exhibit some examples of compact symplectic solvmanifolds $M$
which are $s$--Lefschetz but not $(s+1)$--Lefschetz, for small
values of $s$, and so the map $i$ is bijective for
 $k\geq 2n-s$ and they satisfy the $d\delta$--lemma up to
degree $s$. We detect that they do not satisfy the
$d\delta$--lemma up to degree $s+1$ by exhibiting an invariant
symplectically harmonic $(s+1)$--form $x$ such that $x \in
\Im\delta$ but $x \notin \Im d$. We also find an invariant class
$u\in H^{2n-s-1}_{\delta}(M,\omega)$ such that $i(u)=0$
in~$H^{2n-s-1}(M)$.

\section{Harmonic cohomology of $s$--Lefschetz manifolds} \label{someresults}

We recall some definitions and results about the spaces of
harmonic cohomology classes that we will need in the following
sections. Let $(M,\omega)$ be a symplectic manifold of dimension
$2n$. Denote by $\Omega^*(M)$ the algebra of differential forms on
$M$, by ${\mathfrak X}(M)$ the Lie algebra of vector fields on
$M$, and by ${\cal F}(M)$ the algebra of differentiable functions
on $M$. Since $\omega$ is a non-degenerate $2$--form, we have the
volume form $v_M=\frac {\omega^{n}} {n!}$, and the isomorphism
 $$
 \natural:{\mathfrak X}(M) \longrightarrow \Omega^1(M)
 $$
defined by $\natural(X)=\iota_X(\omega)$ for $X\in {\mathfrak
X}(M)$, where $\iota_X$ denotes the contraction by $X$. We extend
$\natural$ to an isomorphism of graded algebras
$\natural:\bigoplus_{k\geq 0}{\mathfrak X}^k(M) \longrightarrow
\bigoplus_{k\geq 0}\Omega^k(M)$, where ${\mathfrak X}^k(M)$
denotes the space of the skew-symmetric $k$--vectors fields.
Libermann (see \cite{Lib, LM}) defined the {\em symplectic star
operator}
 $$
 *:\Omega^k(M) \longrightarrow \Omega^{2n-k}(M)
 $$
by the condition
 $$
 *(\alpha)=(-1)^{k} \iota_{\natural^{-1}(\alpha)}(v_M).
 $$
This operator can be also defined in terms of the skew-symmetric
bivector field $G$ dual to $\omega$, that is,
$G=-\natural^{-1}(\omega)$. ($G$ is the unique non-degenerate
Poisson structure \cite{Lich} associated with $\omega$.) Denote by
$\Lambda^k(G)$, $k \geq 0$, the associated pairing $\Lambda^k(G):
{\Omega^k(M)\times\Omega^k(M)}\longrightarrow {\cal F}(M)$ which
is $(-1)^k$--symmetric (i.e.\, symmetric for even $k$,
anti-symmetric for odd $k$). Imitating the Hodge star operator for
oriented Riemannian manifolds, Brylinski \cite{Bry} defined the
{\em symplectic star operator} by the condition $\beta\wedge
(*\alpha)=\Lambda^k(G)(\beta, \alpha) v_M$, for $\alpha, \beta \in
\Omega^k(M)$. An easy consequence is that $*^2 = \Id$.

Koszul \cite{Kos} introduced the  differential  $\delta\colon
\Omega^k(M) \longrightarrow \Omega^{k-1}(M)$ on any Poisson
manifold $M$, with Poisson tensor $G$,  by the condition
 $$
 \delta = [\iota_G,d],
 $$
and he proved that $\delta^2=d\delta+\delta d=0$. Later work by
Brylinski \cite{Bry}, shows that the Koszul differential is a {\em
symplectic codifferential} of the exterior differential with
respect to the {\em symplectic star operator}, that is,
 $$
 \delta \alpha =(-1)^{k+1}*d *(\alpha),
 $$
for $\alpha \in \Omega^k(M)$. As in Riemannian Hodge theory, a
$k$--form $\alpha \in \Omega^k(M)$ is said to be {\em coclosed} if
$\delta\alpha=0$, {\em coexact} if $ \alpha=\delta\beta$ for some
$\beta$, and {\em symplectically harmonic} if it is closed and
coclosed. Notice that Koszul definition of $\delta$ implies  that
if $\alpha$ is closed, then $\delta\alpha$ is exact. In \cite{LS}
it is proved the following Leibniz rule for $\delta$. If $f$ is an
arbitrary differentiable function on $M$ and $\alpha \in
\Omega^*(M)$, then
$$
\delta(f\alpha)=f\, \delta \alpha - \iota_{X_{f}}(\alpha),
$$
where $X_f$ is the Hamiltonian vector field of $f$, i.e.,\
$\iota_{X_f} (\omega) = df$.

Let $\Omega^k_{\rm hr}(M,\omega)=\{\alpha\in\Omega^k(M) \mid
d\alpha=\delta\alpha=0\}$ be the space of the symplectically
harmonic $k$--forms. For the de Rham cohomology classes of $M$, we
consider the vector space
 $$
 H^k_{\rm hr}(M,\omega)={\Omega^k_{\rm hr}(M,\omega)\over
 \Omega^k_{\rm hr}(M,\omega)\cap {\rm Im}\, d}\, ,
 $$
consisting of the cohomology classes in $H^k(M)$ containing at
least one symplectically harmonic form.

For $p,k \geq 0$  we define
 $$
 L^p:\Omega^k(M) \longrightarrow \Omega^{2p+k}(M)
 $$
by $L^p(\alpha)=\omega^p \wedge \alpha$ for $\alpha \in
\Omega^k(M)$. In \cite{Ya} it is proved the property following

\begin{lemma}\hskip-.1cm{\rm \cite{Ya}}
(Duality on differential forms). \label{dualitydifforms} The map
 $$
 L^{n-k} \colon \Omega^{k}(M)\longrightarrow
 \Omega^{2n-k}(M)
 $$
is an isomorphism for $0 \leq k\leq n-1$.
\end{lemma}

Since $\omega^p$ is closed, we have
 $$
 [L^p,d] =L^p \circ d - d \circ L^p = 0,
 $$
and the map $L^p$ induces a map $L^p: H^k(M) \longrightarrow
H^{2p+k}(M)$ on cohomology. However, the isomorphisms of Lemma
\ref {dualitydifforms} do not imply special properties on the maps
on cohomology (see Definition \ref{weLefschetz}). Relations
between the operators $\iota_G$, $L$, $d$ and $\delta$ were proved
by Yan in \cite{Ya}. Here we mention the following
 $$
   \iota_{G}=-*L*, \quad
   [\iota_{G},\delta]=0, \quad
   [L,\delta]=-d,
 $$
which implies that if $\alpha$ is coclosed then $d\alpha$ is
coexact, and if $\alpha$ is a symplectically harmonic form then
$L\alpha$ and $\iota_{G}\alpha$ are symplectically harmonic. Also
in \cite{Ya} the following is proved.

\begin{lemma}\hskip-.1cm{\rm \cite{Ya}} (Duality on harmonic forms).
\label{dualityforms} The map
 $$
 L^{n-k} \colon \Omega^{k}_{\rm hr}(M,\omega)\longrightarrow
 \Omega^{2n-k}_{\rm hr}(M,\omega)
 $$
is an isomorphism for $0 \leq k\leq n-1$.
\end{lemma}

Lemma~\ref{dualityforms} implies that the homomorphism
 $$
 L^{n-k}\colon H^{k}_{\rm hr}(M,\omega)\longrightarrow H^{2n-k}_{\rm
 hr}(M,\omega)
 $$
is surjective. (Notice that the duality on harmonic forms may be
not satisfied at the level of the spaces $H^{*}_{\rm
hr}(M,\omega)$.) Since $H^{2n-k}_{\rm hr}(M,\omega)$ is a subspace
of the de Rham cohomology $H^{2n-k}(M)$, we conclude that (see
\cite[Corollary 1.7]{IRTU})
 $$
 H^{2n-k}_{\rm hr}(M,\omega) = \Im(L^{n-k} \colon H^{k}_{\rm
 hr}(M,\omega)\longrightarrow H^{2n-k}(M)).
 $$

A nonzero $k$--form $\alpha$, with $k \leq n$,  is called {\em
primitive} (or {\em effective}) if $L^{n-k+1}(\alpha)=0$. Thus,
any $1$--form is primitive.

\begin{lemma}\hskip-.1cm{\rm \cite[page
46]{LM}} \label{effectiveforms} If $\alpha$ is a primitive
$k$--form, then there is a constant $c$ such that its symplectic
star operator  $*\alpha$ satisfies $*\alpha = c\,
L^{n-k}(\alpha)$.
\end{lemma}

Notice that the previous lemma implies that every closed primitive
$k$--form is symplectically harmonic, and in particular
$H^1(M)=H^{1}_{\rm hr}(M,\omega)$. For the classes in $H^2(M)$,
Mathieu proved that any cohomology class of degree $2$ has a
symplectically harmonic representative, i.e.,\ $H^2(M)=H^{2}_{\rm
hr}(M,\omega)$.

\begin{lemma} \label{primitivedef}
Let $\alpha$ a $k$--form with $k\leq n$.  Then, $\alpha$ is
primitive if and only if $\iota_G(\alpha)=0$.
\end{lemma}

\begin{proof}
It follows from the identity $\iota_G=-*L*$ and Lemma
\ref{effectiveforms}. If $\alpha$  is primitive,
$\iota_G(\alpha)=-*L*(\alpha)=- * c\, L^{n-k+1}(\alpha)=0$.
\end{proof}

\begin{lemma} \label{gradop}
If $\alpha$ is a primitive $k$--form then, for all $j\leq n-k$,
there is a non-zero constant $c_{j,k}$ such that
$\iota_{G}^{j}L^{j}(\alpha) = c_{j,k}\, \alpha$.
\end{lemma}

\begin{proof}
In \cite{Ya} it is proved the relation $[\iota_G,L]=A$, where
$A=\sum (n-k)\pi_k$, $\pi_k$ being the projection. Thus, for $j=1$
we have that $\iota_GL(\alpha)=A\alpha+L(\iota_G\alpha)=
(n-k)\alpha$ because $\alpha$ is primitive. Suppose that
$\iota_{G}^{j}L^{j}(\alpha) = c_{j,k}\, \alpha$ for some $j< n-k$
with $c_{j,k}$ a non-zero constant. Hence,
$\iota_{G}^{j+1}L^{j+1}(\alpha) = \iota_{G}^{j}\iota_G L
L^{j}(\alpha) = \iota_{G}^{j}L\iota_G L^{j}(\alpha)
+(n-k-2j)\iota_{G}^{j}L^{j}(\alpha)= \iota_{G}^{j}L\iota_G
L^{j}(\alpha) +(n-k-2j) c_{j,k}\,(\alpha)$ by the induction
hypothesis. After $p$ times we get that
$$
\iota_{G}^{j+1}L^{j+1}(\alpha) = \iota_{G}^{j-p} L \iota_G^{p+1}
L^{j}(\alpha) + (p+1)(n-k-2j+p)c_{j,k}\,\alpha.
$$
Therefore, for $p=j-1$ and using the induction hypothesis we
conclude that $\iota_{G}^{j+1}L^{j+1}(\alpha) =
c_{j+1,k}\,\alpha$, with $c_{j+1,k}=(j+1)(n-k-j)c_{j,k}$ a
non-zero constant.
\end{proof}

\section{The $d\delta$--lemma for $s$--Lefschetz manifolds} \label{ddeltalemma}

This section is devoted to the study of the $d\delta$--lemma for
symplectic manifolds which are not necessarily hard  Lefschetz. By
Definition \ref{def:ddeltas}, $(M,\omega)$ satisfies the
$d\delta$-lemma up to degree $s$ if $\Im d\cap \ker \delta  = \Im
d \delta =   \Im\delta \cap \ker d$ on $\Omega^k(M)$, for $k \leq
s$ and $\Im d\cap \ker \delta =\Im d \delta$ on $\Omega^{s+1}(M)$.
By applying duality with the symplectic $*$-operator, this is
equivalent to
  \begin{equation}\label{eqn:true2}
 \begin{array}{ll}
 \Im \delta\cap \ker d = \Im d \delta =   \Im d\cap \ker \delta,
  \qquad   & {\mathrm on }\ \Omega^{2n-k}(M),\ {\mathrm for }\ k \leq s, \\
 \Im \delta\cap \ker d =\Im d \delta,
  \qquad   & {\mathrm on }\ \Omega^{2n-s-1}(M).
 \end{array}
  \end{equation}
Let us see the one implication (the other one is proved in an
analogous way). Suppose that $(M,\omega)$ satisfies the
$d\delta$-lemma up to degree $s$. If $\alpha_{2n-k}\in
\Omega^{2n-k}(M)$, $0\leq k \leq s+1$, satisfies that
$\alpha_{2n-k} \in \Im \delta \cap \ker d$, then $*\alpha_{2n-k}$
is a $k$--form in $\Im d \cap\ker \delta=\Im d\delta$, so there is
a $k$--form $\beta_k$ such that $*\alpha_{2n-k}=d\delta (\beta_k)$
and hence $\alpha_{2n-k}=
* d\delta (\beta_k) = - \delta d (*\beta_k) =d\delta (*\beta_k)$.
The equality $\Im d\cap \ker \delta = \Im d \delta$ on
$\Omega^{\geq (2n-s)}(M)$ is proved analogously.

\medskip

Note that if $(M,\omega)$ satisfies the $d\delta$--lemma up to
degree $n-1$ then  both (\ref{eqn:true1}) and (\ref{eqn:true2})
hold for $s=n-1$, and hence $(M,\omega)$ satisfies the
$d\delta$--lemma since then (\ref{eqn:v1a}) also holds on the
space~$\Omega^{n}(M)$.

\medskip

In order to prove Theorems \ref{gradosddelta} and \ref{mapi} we
need the following results.

\begin{lemma} \label{ddeltayL}
Let $(M,\omega)$ be a symplectic manifold of dimension $2n$, and
let $\alpha$ be a $k$--form. Then
 \begin{enumerate}
 \item[{\rm (i)}] $d\delta(L^{p}(\alpha))=L^{p}(d\delta(\alpha))$ for all $p\geq 0$.
 \item[{\rm (ii)}] If $\alpha$ is primitive, then $d \delta
 (\alpha)$ is also primitive.
 \end{enumerate}
\end{lemma}

\begin{proof}
Since $[L,\delta]=-d$, we see that $\delta L=L \delta+d$. Thus,
$d\delta(L^{p}(\alpha))=d(L \delta+d)L^{p-1}(\alpha)=dL \delta
L^{p-1}(\alpha)$. Proceeding in this fashion $p$ times, and using
that $L$ and $d$ commute, we obtain (i). Now to show (ii) we have,
using (i), that $L^{n-k+1}(d\delta(\alpha))
=d\delta(L^{n-k+1}(\alpha)) =0$ since $\alpha$ is a primitive
$k$--form.
\end{proof}

\begin{lemma} \label{deltaexact}
Let $(M,\omega)$ be a symplectic manifold of dimension $2n$, let
$\beta$ be a $r$--form, and let $\alpha=L^{p}(\beta)$, with $p\geq
0$. If $\delta \beta$ is exact, then $\delta \alpha$ is also
exact.
\end{lemma}

\begin{proof}
Write $\delta\beta=d\gamma$. Using $[L,\delta]=-d$, we have
$\delta \alpha =\delta L^{p}( \beta)=(L\delta+d)L^{p-1}( \beta)=
L\delta L^{p-1}( \beta)+dL^{p-1}( \beta)$. Proceeding in a similar
way with the first summand, after $p$ steps,  we get
 $$
 \delta  \alpha=\delta L^{p}( \beta)=L^{p}\delta ( \beta)+
 p\, d L^{p-1}( \beta)= d(L^{p}(\gamma)+pL^{p-1}(\beta)),
 $$
which proves the lemma.
\end{proof}

Consider  a $(2n-i)$--form $\psi$ on $(M,\omega)$ with $i\leq n$.
According to the duality on differential forms,  there is a unique
$i$--form $\varphi$ such that $\psi = L^{n-i}(\varphi)$. Lepage
decomposition theorem \cite{LM} implies that $\varphi$ may be
uniquely decomposed as a sum
 \begin{equation} \label{eqn:10}
 \varphi = \varphi_{i} +L(\varphi_{i-2})+\cdots + L^{q}(\varphi_{i-2q}),
 \end{equation}
with $q\leq [i/2]$, where $[i/2]$ being the largest integer less
than or equal to $i/2$, and where the form $\varphi_{i-2j}$ is a
primitive $(i-2j)$--form, for $j=0,\ldots,q$. This implies that
$\psi= L^{n-i}(\varphi)$ may be uniquely decomposed as  the sum
 \begin{equation} \label{eqn:10a}
 \psi = L^{n-i}(\varphi)= L^{n-i}(\varphi_{i})
 +L^{n-i+1}(\varphi_{i-2})+\cdots + L^{n-i+q}(\varphi_{i-2q}).
 \end{equation}

\begin{lemma} \label{gradoalto}
Let $(M,\omega)$ be a symplectic manifold of dimension $2n$, and
let  $\psi=L^{n-i}(\varphi)\in \Omega^{2n-i}(M)$ with $i\leq n$.
 \begin{enumerate}
 \item[{\rm (i)}] If  $d \delta (\psi)=0$,
or equivalently $d \delta (\varphi)=0$, then all the forms
$\varphi_{i-2j}$ in the decomposition {\rm (\ref {eqn:10})} and
 {\rm (\ref {eqn:10a})} satisfy $d \delta (\varphi_{i-2j})=0$.
 \item[{\rm (ii)}]  If $\delta \varphi_{i-2j}$
is exact for all $j=0,\ldots,q$, then both $\delta \varphi$ and
$\delta \psi$ are exact.
 \end{enumerate}
\end{lemma}

\begin{proof}
Suppose that $d \delta (\psi)=0$. Applying $d\delta$ to
(\ref{eqn:10a}), using Lemma \ref{ddeltayL} and the uniqueness of
the decomposition, we have that
 $$
 L^{n-i+j} d \delta (\varphi_{i-2j})=0,
 $$
for $j=0,\ldots,q$. We see that $L^{n-i+j} d \delta
(\varphi_{i-2j})=0$ implies $d \delta (\varphi_{i-2j})=0$. In
fact, the map $L^{n-i+2j}\colon \Omega^{i-2j}(M)\longrightarrow
\Omega^{2n-i+2j}(M)$ is an isomorphism for all $j=0,\ldots,q$. So,
the map $L^{n-i+j}\colon \Omega^{i-2j}(M)\longrightarrow
\Omega^{2n-i}(M)$ is injective for $j=1,\ldots,q$, and it is an
isomorphism for $j=0$. Using again Lemma \ref{ddeltayL} and the
duality on differential forms, one can check that $d \delta
(\varphi)=0$ implies the same result. Part (ii) follows from Lemma
\ref{deltaexact} and using that $\delta$ is a linear map.
\end{proof}

\begin{proposition} \label{s-Lefschetzandddelta}
Let $(M,\omega)$ be an $s$--Lefschetz compact symplectic manifold
of dimension $2n$.  Then,
 $$
 \Im\delta\cap \ker d  = \Im d\cap  \Im\delta,
 $$
 on the spaces $\Omega^{\leq s}(M)$ and $\Omega^{\geq (2n-s-2)}(M)$; and
 $$
 \Im d\cap\ker \delta=\Im d\cap  \Im\delta,
 $$
on $\Omega^{\leq (s+2)}(M)$ and $\Omega^{\geq (2n-s)}(M)$.
\end{proposition}

\begin{proof}
We prove only the first identity because the second is analogous
by duality using the symplectic $*$-operator. The result can be
restated in the following way: if $\varphi$ is a $k$--form, with
$k\leq s+1$ or $k\geq 2n-s-1$, and such that $d\delta(\varphi)=0$,
then $\delta\varphi$ is exact.

First, we  show such  a result for any  primitive $k$--form
$\varphi$ with $k\leq s+1$. We define the $(k-1)$--form $\gamma$
by
 \begin{equation}\label{eqn:01}
 L^{n-k+1}(\gamma)=d L^{n-k}(\varphi).
 \end{equation}
Thus $\gamma$ is primitive since $L^{n-k+2}(\gamma)=d
L^{n-k+1}(\varphi) =0$. Applying $\iota_{G}^{n-k+1}$ in
(\ref{eqn:01}), using Lemma~\ref{gradop} and
$\delta=[\iota_{G},d]$, we have
 $$
 c_{n-k+1,k-1}\gamma =\iota_{G}^{n-k+1}dL^{n-k}\varphi=
 \iota_{G}^{n-k}(d\iota_{G}+\delta)L^{n-k} \varphi.
 $$
Proceeding in this fashion, after $(n-k+1)$ times, we have
 $$
 c_{n-k+1,k-1}\gamma = (d\iota_{G}^{n-k+1} -
 (n-k+1)\delta\iota_{G}^{n-k})L^{n-k}\varphi.
 $$
Since $\varphi$ is primitive, $(d\iota_{G}^{n-k+1})L^{n-k}\varphi
= d\iota_{G}(c_{n-k,k}\varphi)=0$ by Lemma \ref{primitivedef}. So,
there is a non-zero constant $c$ such that $\gamma= c\,  \delta
\varphi$. Applying $L^{n-k+1}$ to both sides and using
(\ref{eqn:01}) we obtain
 $$
 c\, L^{n-k+1}\delta\varphi = L^{n-k+1}\gamma=d(L^{n-k}\varphi).
 $$
By hypothesis $\delta\varphi$ is closed. Moreover the map
$L^{n-k+1}\colon H^{k-1}(M)\longrightarrow H^{2n-k+1}(M)$ is an
isomorphism for $k-1\leq s$ since $(M,\omega)$ is compact and
$s$--Lefschetz. Thus  $\delta\varphi$ is exact because
$L^{n-k+1}\delta\varphi$ defines the zero class.

Now we pass to the case where $\varphi$ is an arbitrary $k$--form
with $k\leq s+1$ such that $d \delta(\varphi)=0$. {}From Lemma
\ref{gradoalto} we know that every primitive form $\varphi_{i-2j}$
in the decomposition (\ref{eqn:10}) satisfies $d
\delta(\varphi_{i-2j})=0$, and so $\delta(\varphi_{i-2j})$ is
exact. Now  Lemma \ref{gradoalto} implies that $\delta \varphi$ is
exact.

Finally, if $\psi$ is a $k$--form with $k\geq 2n-s-1$ and such
that $d \delta(\psi)=0$, then the forms $\varphi_{i-2j}$ in the
decomposition (\ref{eqn:10a}) are of degree $\leq s+1$, and they
satisfy $d\delta(\varphi_{i-2j})=0$ by Lemma \ref{gradoalto}.
Taking account the previous result for primitive forms, we
conclude that all the forms $\delta \varphi_{i-2j}$ are exact, and
hence $\delta \psi$ is exact by Lemma \ref{gradoalto}.
\end{proof}

\begin{proposition} \label{s-Lefschetzandddelta1}
Let $(M,\omega)$ be an $s$--Lefschetz compact symplectic manifold
of dimension $2n$.  We have
 \begin{enumerate}
 \item[{\rm (i)}] $\Im \delta \cap \ker d = \Im d \cap \ker\delta$ on
$\Omega^{\leq s}(M)$ and $\Omega^{\geq (2n-s)}(M)$.
 \item[{\rm (ii)}] $(M,\omega)$ satisfies the $d\delta$--lemma
 up to  degree $s$.
  \end{enumerate}
\end{proposition}

\begin{proof}
Part (i) follows directly from Proposition
\ref{s-Lefschetzandddelta}.

To show (ii), we shall first prove that $\Im \delta \cap \ker
d=\Im d \delta$ on the spaces $\Omega^{\geq (2n-s-1)}(M)$. We will
prove this by induction on $s$. For $s=0$, assume $\alpha \in
\Omega^{2n}(M)$ such that $\alpha \in \Im \delta\cap\ker d   =
 \Im d \cap\ker \delta$. Then $\alpha=0$ because $\Im \delta=0$,
and so $\alpha=0=d\delta0$. Now we see that $\Im\delta\cap\ker d =
\Im d \delta$ on $\Omega^{2n-1}(M)$. Let $\beta=\delta\alpha$ be a
$(2n-1)$--form, $\alpha \in \Omega^{2n}(M)$, such that
$d\delta\alpha=0$. Since $d\alpha=0$ and $H^{2n}(M)=H^{2n}_{\rm
hr}(M,\omega)$, there is $\tilde\alpha\in\Omega^{2n}(M)$ such that
$\delta\tilde\alpha=0$, $d\tilde\alpha=0$ and
$\alpha=\tilde\alpha+d\gamma$, for some
$\gamma\in\Omega^{2n-1}(M)$. Therefore, $\beta=\delta\alpha=\delta
d\gamma=d\delta(-\gamma)$.

Now take $s>0$, and assume that if $(M,\omega)$ is
$(s-1)$--Lefschetz, then $\Im \delta \cap \ker d  = \Im d \delta$
on $\Omega^{\geq (2n-s)}(M)$. We need to prove that if
$(M,\omega)$ is $s$--Lefschetz, then $\Im \delta \cap \ker d  =
\Im d \delta$ on $\Omega^{2n-s-1}(M)$. We will use subscripts to
keep track of the spaces that the forms belong to, i.e.\
$\alpha_{k} \in \Omega^{k}(M)$. We consider a $(2n-s-1)$--form
$\alpha_{2n-s-1}$ such that $\alpha_{2n-s-1}= \delta \alpha_{2n-s}
\in \Im \delta \cap \ker d$. Then,
  $$
  0 = d \alpha_{2n-s-1} = d\delta \alpha_{2n-s} = -\delta d\alpha_{2n-s},
  $$
which implies that $d\alpha_{2n-s}$ is a $(2n-s+1)$--form such
that $d\alpha_{2n-s} \in \Im d \cap \ker \delta = \Im \delta \cap
\ker d = \Im d \delta$ by (i) and induction hypothesis. Thus
 $$
 d\alpha_{2n-s} = d\delta \alpha_{2n-s+1},
 $$
for some $\alpha_{2n-s+1} \in \Omega^{2n-s+1}(M)$, and
consequently
 $$
 d(\alpha_{2n-s}-\delta \alpha_{2n-s+1})=0,
 $$
which means that   $(\alpha_{2n-s} - \delta \alpha_{2n-s+1})$
defines a de Rham cohomology class in $H^{2n-s}(M) = H^{2n-s}_{\rm
hr}(M,\omega)$, the last equality by Theorem~\ref{grados+2}. Thus,
there exist a symplectically harmonic $(2n-s)$--form
$\beta_{2n-s}$ and $\eta_{2n-s-1} \in \Omega^{2n-s-1}(M)$  such
that
 $$
 \alpha_{2n-s} - \delta \alpha_{2n-s+1} - \beta_{2n-s} = d\eta_{2n-s-1}.
 $$
Applying $\delta$ to both sides we have
 $$
 \alpha_{2n-s-1} = \delta \alpha_{2n-2}= \delta d\eta_{2n-s-1}=-d\delta \eta_{2n-s-1}
 \in \Im d\delta.
 $$

To end the proof, we use the duality by the symplectic
$*$-operator to show that $\Im d\cap \ker \delta =\Im d\delta$ on
the spaces $\Omega^{\leq (s+1)}(M)$. In fact, let us consider
$\alpha_r$ a differential $r$--form, with $r\leq s+1$, such that
$\alpha_r \in \Im d\cap \ker \delta $. Then, $*\alpha_r$ is a
$(2n-r)$--form, $2n-r \geq 2n-s-1$, such that $*\alpha_r \in \Im
\delta\cap \ker d  = \Im d \delta$, and so $\alpha_r \in \Im d
\delta$. The equality $\Im d\cap \ker \delta =\Im \delta\cap \ker
d$ on the spaces  $\Omega^{\leq s}(M)$ follows from (i), and this
completes the proof of the $d\delta$-lemma up to degree $s$ for
$(M,\omega)$.
 \end{proof}

\bigskip

\noindent{\bf Proof of Theorem~\ref{gradosddelta}~:} Clearly (i)
implies (ii) by Proposition \ref{s-Lefschetzandddelta1}. Also (ii)
implies (iii) by duality of the symplectic $*$-operator.

Let us show that (iii) implies (i). By Theorem \ref{grados+2}, it
is enough to prove that every de Rham cohomology class of degree
$k$ has a symplectically harmonic representative for $2n-s\leq
k\leq 2n$. Let us consider $[\gamma] \in H^{k}(M)$ with $2n-s\leq
k\leq 2n$. Then $d\gamma=0$, and $\delta \gamma$ is a
$(k-1)$--form such that $d\delta \gamma=0$ since $d$ and $\delta$
anticommute. This means that $\delta \gamma$ lives in
$\Im\delta\cap\ker d$ which is equal to $\Im d\delta$ on forms of
degree $k-1 \geq 2n-s-1$ by the hypothesis (iii). This implies
that there is a $(k-1)$--form $\theta$ such that $\delta
\gamma=d\delta\theta$. So $\delta(\gamma+d\theta)=0$. Then, the
form $\gamma+d\theta$ is symplectically harmonic and cohomologous
to $\gamma$.\QED

\begin{remark}
Notice that if $(M,\omega)$ is a compact symplectic manifold of
dimension $2n$ and it is $(n-2)$--Lefschetz, then the identities
{\rm (\ref{eqn:v1a})} hold on $\Omega^{\leq (n-2)}(M)$ and
$\Omega^{\geq (n+2)}(M)$, and also $\Im \delta\cap \ker d=\Im
d\cap \Im\delta=\Im d\cap \ker \delta$ on $\Omega^n(M)$, by
Proposition \ref{s-Lefschetzandddelta}. Nonetheless, if
$(M,\omega)$ is not hard Lefschetz, then this last space is in
general different from \, $\Im d\delta$.
\end{remark}

\bigskip

Let $\Omega^k_{\delta}(M, \omega)=\{\alpha\in\Omega^k(M) \mid
\delta\alpha=0\}$ be the space of the coclosed $k$--forms. Since
$d$ and $\delta$ anti-commute, then
$d(\Omega^k_{\delta}(M,\omega))\subset
\Omega^{k+1}_{\delta}(M,\omega)$, and so
$(\Omega^*_{\delta}(M,\omega),d)$ is a subcomplex of the de Rham
complex $(\Omega^*(M),d)$. We denote by $H^*_{\delta}(M,\omega)$
its cohomology, that is
 $$
 H^k_{\delta}(M,\omega)={\ker(d \colon \Omega^k_{\delta}(M,\omega)
 \longrightarrow \Omega^{k+1}_{\delta}(M,\omega))
 \over
 \Im(d \colon \Omega^{k-1}_{\delta}(M,\omega)
 \longrightarrow \Omega^{k}_{\delta}(M,\omega))}.
 $$
 Therefore, any cohomology class on $H^k_{\delta}(M,\omega)$
 is symplectically harmonic, and we have a natural map
 $ i_{1} \colon H^k_{\delta}(M,\omega)\longrightarrow H^{k}_{\rm hr}(M,\omega)$
 which is always surjective but may be non-injective.
The next theorem gives a necessary and sufficient
 condition for the injectivity of the map
 $i_{1}$.
(Notice that $\Omega^*_{\delta}(M,\omega)= \bigoplus \,
\Omega^{k}_{\delta}(M,\omega)$ is a vector space but not an
algebra because the codifferential $\delta$ does not satisfy a
Leibniz rule.) It is clear that there is a natural map
 $$
 i \colon H^k_{\delta}(M,\omega)\longrightarrow H^{k}(M),
 $$
for all $k$. In fact, denote by $i_{2}$ the natural
 inclusion
 $$
 i_{2} \colon H^k_{\rm hr}(M,\omega)\longrightarrow H^{k}(M).
 $$
 Then, $i=i_{2}\circ i_{1}$.

 \bigskip

\noindent{\bf Proof of Theorem~\ref{mapi}~:} Suppose that
 $(M,\omega)$ is $s$--Lefschetz. By Theorem \ref{grados+2},
 $H^{k}_{\rm hr}(M,\omega)=H^{k}(M)$ for $k \leq s+2$
 and $k \geq 2n-s$. Then, to show (ii) it is enough to
 prove that the map
 $i=i_{1} \colon H^k_{\delta}(M,\omega)\longrightarrow H^{k}_{\rm hr}(M,\omega)$
 is injective for $k \leq s+1$
 and $k \geq 2n-s$ because such a map is always surjective.
 Consider $[\alpha] \in H^k_{\delta}(M,\omega)$ and suppose
 that $[\alpha]=i[\alpha]$ defines the zero class  on $H^{k}_{\rm hr}(M,\omega)$.
 Then $\alpha$ is exact, i.e.\ $\alpha = d\beta$ for
 some $\beta \in \Omega^{k-1}(M)$. But if $k \leq s+1$
or $k \geq 2n-s$, Theorem \ref{gradosddelta} implies
$\alpha=d\delta\eta$ for some $\eta \in \Omega^{k}(M)$. Hence
$\alpha=d\delta\eta \in \Im (d \colon
\Omega^{k-1}_{\delta}(M,\omega) \longrightarrow
\Omega^{k}_{\delta}(M,\omega))$. This means that $\alpha$ defines
the zero class on $H^k_{\delta}(M,\omega)$, which proves (ii).

Clearly (ii) implies (iii). We show that (iii) implies (i). In
fact, if $[\alpha]\in H^{k}_{\delta}(M,\omega)$, $[\alpha]$ is a
harmonic cohomology class. Thus, if the map $ i \colon
H^k_{\delta}(M,\omega)\longrightarrow H^{k}(M)$ is bijective for
$k \geq 2n-s$ then $H^{k}_{\rm hr}(M,\omega)=H^{k}(M)$ for $k \geq
2n-s$, i.e.\ $(M,\omega)$ is $s$--Lefschetz according to
Theorem~\ref{grados+2}.\QED

\section{Harmonic cohomology of compact
completely solvmanifolds} \label{solvmanifolds}

Let ${\frak g}$ be a Lie algebra of dimension~$2n$, and denote by
$d$ the Chevalley-Eilenberg differential of~${\frak g}$. An
element $\omega\in \bigwedge^2({\frak g^{*}})$ such that
$d\omega=0$ and $\omega^n\not=0$ will be called a {\it symplectic
form} on~${\frak g}$.

Symplectic Hodge theory can be introduced for a symplectic form
$\omega$ on a Lie algebra ${\frak g}$ in a similar way as in
Section~\ref{someresults}. Let us define the star operator $*
\colon \bigwedge^k({\frak g^{*}}) \longrightarrow \bigwedge^{2n-k}
({\frak g^{*}})$ by
 $$
 *\, \alpha=(-1)^{k} \iota_{\natural^{-1}(\alpha)}\frac {\omega^{n}}
 {n!},
 $$
for any $\alpha \in \bigwedge^k({\frak g^{*}})$, where $\natural$
denotes the isomorphism between $\bigwedge^k({\frak g})$ and
$\bigwedge^k({\frak g^{*}})$ extended from the natural isomorphism
$\natural\colon {\frak g} \longrightarrow {\frak g}^*$ given by
$\natural(X)(Y)=\omega(X,Y)$, for $X,Y\in {\frak g}$.

We define the codifferential $\delta\colon \bigwedge^{k} ({\frak
g^{*}}) \longrightarrow \bigwedge^{k-1} ({\frak g^{*}})$ by
$$
 \delta\alpha =(-1)^{k+1}*d *\alpha,
$$
for any $\alpha \in \bigwedge^k ({\frak g^{*}})$. Now, let
$\bigwedge^k_{\rm hr}({\frak g^{*}},\omega)=\{\alpha\in
\bigwedge^k({\frak g^{*}}) \mid d\alpha=\delta\alpha=0\}$, and
consider the space
 $$
 H^k_{\rm hr}({\frak g},\omega)={\bigwedge\phantom{i}\!\!\!^k_{\rm hr}({\frak g^{*}},\omega)\over
 \bigwedge\phantom{i}\!\!\!^k_{\rm hr}({\frak g^{*}},\omega)\cap {\rm Im}\,
 d}\, .
 $$
Then, $H^k_{\rm hr}({\frak g},\omega)$ consists of all the classes
in the Chevalley-Eilenberg cohomology $H^k({\frak g})$ of ${\frak
g}$ containing at least one representative which is both closed
and $\omega$-coclosed.

Let $G\in \bigwedge^2({\frak g})$ be given by
$G=-\natural^{-1}(\omega)$. In order to study the spaces $H^k_{\rm
hr}({\frak g},\omega)$ we consider the linear maps $L\colon
\bigwedge^*({\frak g^{*}})\longrightarrow \bigwedge^{*+2}({\frak
g^{*}})$, $\iota_G\colon \bigwedge^*({\frak g^{*}})\longrightarrow
\bigwedge^{*-2}({\frak g^{*}})$ and $A\colon \bigwedge^*({\frak
g^{*}})\longrightarrow \bigwedge^*({\frak g^{*}})$ as usual:
$L\alpha$ is the wedge product by $\omega$, $\iota_G\alpha$ the
contraction by $G$ and $A=\sum(n-k)\pi_k$, where $\pi_k$ is the
projection onto $\bigwedge^{k}({\frak g^{*}})$.
Following~\cite{Ya}, although the arguments in this special case
are more direct, the following relations hold:
 $$
 [L,\delta]=-d,\quad [\iota_G,d]=\delta,\quad
 [L,d]=[\iota_G,\delta]=0,
 $$
and
 $$
 [\iota_G,L]=A, \quad [A,\iota_G]=2\, \iota_G, \quad [A,L]=-2\,L.
 $$

Since the standard basis $\left\{X=\left(\begin{array}{cc}
0&1\\
0&0
\end{array}\right),\ Y=\left(\begin{array}{cc}
0&0\\
1&0
\end{array}\right),\ H=\left(\begin{array}{cc}
1&0\\
0&-1
\end{array}\right)\right\}$ of ${\frak s}{\frak l}(2,\CC)$ satisfies
 $$
 [X,Y]=H,\quad [H,X]=2X,\quad [H,Y]=-2Y,
 $$
we have representations $\rho_1\colon {\frak s}{\frak l}(2,\CC)
\longrightarrow {\frak g}{\frak l}(\bigwedge^*({\frak
g}^*)\otimes\CC)$ and $\rho_2\colon {\frak s}{\frak l}(2,\CC)
\longrightarrow {\frak g}{\frak l}(\bigwedge^*_{\rm hr}({\frak
g}^*,\omega)\otimes\CC)$ of the Lie algebra ${\frak s}{\frak
l}(2,\CC)$ on the complex vector spaces $\bigwedge^*({\frak
g}^*)\otimes\CC$ and $\bigwedge^*_{\rm hr}({\frak
g}^*,\omega)\otimes\CC$, respectively, defined by the
correspondence
 $$
 \rho_i(X)=\iota_G,\quad \rho_i(Y)= L,\quad \rho_i(H)= A \quad
 (i=1,2),
 $$
where $\iota_G$, $L$ and $A$ are understood for $\rho_1$ as the
extension of the maps $\iota_G$, $L$ and $A$ above to the
complexification $\bigwedge^*({\frak g}^*)\otimes\CC$ of
$\bigwedge^*({\frak g}^*)$, and for $\rho_2$ as the restriction of
them to the subspace $\bigwedge^*_{\rm hr}({\frak
g}^*,\omega)\otimes\CC$. Notice that we can consider the
restriction $\rho_2$ of the ${\frak s}{\frak l}(2,\CC)$
representation $\rho_1$ since if $\alpha$ is symplectically
harmonic then $L\alpha$ and $\iota_G\alpha$ are symplectically
harmonic.

It is well-known (see for example~\cite{wells}) that for any
representation $\rho$ of ${\frak s}{\frak l}(2,\CC)$ on a finite
dimensional complex vector space $V$, all the eigenvalues of
$\rho(H)\colon V\longrightarrow V$ are integer numbers and, if
$V_k$ denotes the eigenspace of $\rho(H)$ with respect to the
eigenvalue $k$, then
 $$
 \rho(Y)^k\colon V_{-k} \longrightarrow V_{k} \quad \mbox{ and }
 \quad \rho(X)^k\colon V_{k} \longrightarrow V_{-k}
 $$
are isomorphisms. Therefore, since $\bigwedge^r({\frak
g^{*}})\otimes\CC$ and $\bigwedge^r_{\rm hr}({\frak
g^{*}},\omega)\otimes\CC$ are the eigenspaces of $\rho_1(H)$ and
$\rho_2(H)$, respectively, with respect to the eigenvalue $r$, we
conclude that
 $$
 L^k\colon \bigwedge\phantom{i}\!\!\!^{n-k}({\frak g^{*}})
 \longrightarrow \bigwedge\phantom{i}\!\!\!^{n+k}({\frak g^{*}})
 $$
and
 $$
 L^k\colon \bigwedge\phantom{i}\!\!\!^{n-k}_{\rm hr}({\frak
 g^{*}},\omega) \longrightarrow
 \bigwedge\phantom{i}\!\!\!^{n+k}_{\rm hr}({\frak g^{*}},\omega)
 $$
are isomorphisms for $k\geq 0$.

\begin{remark}\label{duality-idea}
Lemmas~\ref{dualitydifforms} and~\ref{dualityforms} expressing
duality of forms and of harmonic forms, respectively, are derived
by Yan~\cite{Ya} from the theory of a special type of {\it
infinite dimensional} representations of ${\frak s}{\frak
l}(2,\CC)$ called {\it of finite $H$-spectrum}. Any finite
dimensional representation of ${\frak s}{\frak l}(2,\CC)$ is of
this type.
\end{remark}

The following result is a direct consequence of the isomorphisms
$L^k$ given above. In the proof we follow the lines of~\cite[Lemma
4.3]{Yam} and~\cite[Corollary 2.4]{IRTU}, where a similar result
is given for the harmonic cohomology $H^{*}_{\rm hr}(M)$ of a
symplectic manifold $M$.

\begin{lemma}\label{desc-harm-g}
Let ${\frak g}$ be a $2n$-dimensional Lie algebra with a
symplectic form $\omega$. For every $k\ge 0$, we have
 $$
  H^{n-k}_{\rm hr}({\frak g},\omega) = P_{n-k}({\frak
 g},\omega)+L(H^{n-k-2}_{\rm hr}({\frak g},\omega)) \quad\mbox{ and
 }\quad H^{n+k}_{\rm hr}({\frak g},\omega) = L^k(H^{n-k}_{\rm
 hr}({\frak g},\omega)),
 $$
where $P_r({\frak g},\omega)=\{[\alpha]\in H^r({\frak g}) \mid
L^{n-r+1}[\alpha]=0 \}$ is the space of primitive cohomology
classes of degree $r$, and $L$ denotes the product by $[\omega]\in
H^2({\frak g})$.
\end{lemma}

\begin{proof}
Let $[\alpha]\in H^{n-k}_{\rm hr}({\frak g},\omega)$. Since
$L^{k+2}(\bigwedge^{n-k-2}_{\rm hr}({\frak
g^{*}},\omega))=\bigwedge^{n+k+2}_{\rm hr}({\frak g^{*}},\omega)$,
there exists $\beta$ such that $d\beta=\delta\beta=0$ and
$L^{k+2}\beta=L^{k+1}\alpha$. Therefore,
$L^{k+1}([\alpha]-L[\beta])=0$. Since $[\alpha] =
([\alpha]-L[\beta]) + L[\beta]$, the inclusion $H^{n-k}_{\rm
hr}({\frak g},\omega) \subset P_{n-k}({\frak
g},\omega)+L(H^{n-k-2}_{\rm hr}({\frak g},\omega))$ holds.

To prove the other inclusion it suffices to show that any class
$[\alpha]\in P_{n-k}({\frak g},\omega)$ contains a representative
$\tilde\alpha$ such that $\delta\tilde\alpha=0$. Since
$L^{k+1}[\alpha]=0$, there exists $\gamma\in
\bigwedge^{n+k+1}({\frak g}^*)$ such that $L^{k+1}\alpha=d\gamma$,
so $\gamma=L^{k+1}\beta$ for some $\beta\in
\bigwedge^{n-k-1}({\frak g}^*)$. Let $\tilde\alpha=\alpha-d\beta$.
Since $L^{k+1}\tilde\alpha=0$ we have that $*\tilde\alpha$ is
proportional to $L^{k}\tilde\alpha$, therefore $\tilde\alpha$ is a
representative of $[\alpha]$ satisfying
$\delta\tilde\alpha=[\iota_G,d]\tilde\alpha=0$.

Finally, if $[\alpha]\in H^{n+k}_{\rm hr}({\frak g},\omega)$ then
there is $\beta\in\bigwedge^{n-k}_{\rm hr}({\frak g^{*}},\omega)$
such that $\alpha=L^k\beta$, so $H^{n+k}_{\rm hr}({\frak
g},\omega) = L^k(H^{n-k}_{\rm hr}({\frak g},\omega))$.
\end{proof}

Suppose that a simply connected Lie group $G$ has a discrete
subgroup $\G$ such that the quotient $M=\G\backslash G$ is
compact. Let us denote by ${\frak g}$ the Lie algebra of $G$.
Since any element in $\bigwedge^k({\frak g^{*}})$ is identified to
a left invariant form on $G$, it descends to the quotient $M$ and
there is a natural injection $\bigwedge^*({\frak
g^{*}})\hookrightarrow \Omega^*(M)$ which commutes with the
differentials.

On the other hand, if the Lie algebra ${\frak g}$ of $G$ possesses
a symplectic form $\omega$ then it descends to a symplectic form
on $M$, which we shall also denote by $\omega$. In this case the
natural injection $\bigwedge^*({\frak g^{*}})\hookrightarrow
\Omega^*(M)$ also commutes with the symplectic stars, and so with
the $\delta$'s. Therefore, we have a natural homomorphism
$H^*_{\rm hr}({\frak g},\omega)\longrightarrow H^*_{\rm hr}(M)$.

\begin{proposition}\label{general-iso}
If the natural inclusion $\bigwedge^*({\frak
g^{*}})\hookrightarrow \Omega^*(M)$ induces an isomorphism
$H^*({\frak g})\cong H^*(M)$ in cohomology, then the inclusion
$\bigwedge^*_{\rm hr}({\frak g^{*}},\omega)\hookrightarrow
\Omega^*_{\rm hr}(M)$ also induces an isomorphism $H^*_{\rm
hr}({\frak g},\omega)\cong H^*_{\rm hr}(M)$.
\end{proposition}

\begin{proof}
Since the natural homomorphism $H^k({\frak g})\longrightarrow
H^k(M)$ commutes with $L$ and it is an isomorphism, for each
$k\leq n$ we have an isomorphism between $P_k({\frak g},\omega)$
and the space $P_k(M)=\{ [\alpha]\in H^k(M) \mid
L^{n-k+1}[\alpha]=0 \}$ of primitive cohomology classes of degree
$k$. Since $H^{n+k}_{\rm hr}(M) = L^k(H^{n-k}_{\rm hr}(M))$, from
Lemma~\ref{desc-harm-g} it suffices to prove that $H^{n-k}_{\rm
hr}({\frak g},\omega) \cong H^{n-k}_{\rm hr}(M)$. But this follows
easily by an inductive argument, taking into account
Lemma~\ref{desc-harm-g} and the fact~\cite{Yam} that $H^{n-k}_{\rm
hr}(M) = P_{n-k}(M)+L(H^{n-k-2}_{\rm hr}(M))$. Notice that for
starting the induction, i.e.\ for $n-k=0,1,2$, we have
$H^{n-k}_{\rm hr}(M)=H^{n-k}(M)\cong H^{n-k}({\frak
g})=H^{n-k}_{\rm hr}({\frak g},\omega)$.
\end{proof}

Let $M=\G\backslash G$ be a compact solvmanifold, that is, a
compact quotient of a simply connected solvable Lie group $G$ by a
discrete subgroup $\G$. Suppose in addition that the Lie algebra
${\frak g}$ of $G$ is completely solvable, i.e.\ ${\rm ad}_X\colon
{\frak g} \longrightarrow {\frak g}$ has only real eigenvalues for
any $X \in {\frak g}$. By Hattori theorem~\cite{Hat}, which is a
generalization to the completely solvable context of Nomizu
theorem~\cite{No} for nilmanifolds, the natural inclusion
$\bigwedge^*({\frak g^{*}})\hookrightarrow \Omega^*(M)$ induces an
isomorphism $H^*({\frak g})\cong H^*(M)$.

Let $\omega$ be a symplectic form on $M=\G\backslash G$. {}From
the results above it is clear that the harmonic cohomology only
depends on the cohomology class $[\omega]$ of the symplectic form.
Since $H^2(M)\cong H^2({\frak g})$ we can suppose without loss of
generality that $\omega$ is invariant, that is, it stems from a
symplectic form on the Lie algebra ${\frak g}$.

\begin{corollary}\label{Nomizu-type-th}
Let $M=\G\backslash G$ be a compact solvmanifold endowed with a
symplectic form~$\omega$. If the Lie algebra ${\frak g}$ of $G$ is
completely solvable, then the natural injection $\bigwedge^*_{\rm
hr}({\frak g^{*}},\omega)\hookrightarrow \Omega^*_{\rm hr}(M)$
induces an isomorphism $H^*_{\rm hr}({\frak g},\omega)\cong
H^*_{\rm hr}(M)$.
\end{corollary}

In particular, the result holds for symplectic nilmanifolds, which
has been already obtained in~\cite{Yam}.

{}From Theorems~\ref{gradosddelta}  and \ref{mapi} and their
corresponding analogues for a Lie algebra endowed with a
symplectic form, we have the following results.

\begin{corollary}\label{d-delta-invariant}
Let $M=\G\backslash G$ be a compact solvmanifold endowed with a
symplectic form~$\omega$. Then, the $d\delta$-lemma holds on $M$
up to degree~$s$ if and only if it holds on ${\frak g}$ up to
degree~$s$, i.e.,
  $$
 \begin{array}{ll}
 d(\bigwedge\phantom{i}\!\!\!^{k-1}({\frak g}^*)) \cap \ker \delta
 =   \delta (\bigwedge\phantom{i}\!\!\!^{k+1}({\frak g}^*))\cap \ker d=
 d \delta (\bigwedge\phantom{i}\!\!\!^k({\frak g}^*)),
  \qquad   & {\mathrm for }\ k \leq s, \\
 d(\bigwedge\phantom{i}\!\!\!^{s}({\frak g}^*)) \cap \ker \delta =
 d \delta (\bigwedge\phantom{i}\!\!\!^{s+1}({\frak g}^*)).
 \end{array}
  $$
\end{corollary}

\begin{corollary}\label{mapi-invariant}
Let $M=\G\backslash G$ be a compact solvmanifold endowed with a
symplectic form~$\omega$. Then, the map $i$ given in {\rm
(\ref{eqn:v1aa})} is bijective for all $k\geq 2n-s$ if and only if
the map $i\colon H^k_{\delta}({\frak g},\omega)\longrightarrow
H^{k}({\frak g})$ is bijective for all $k\geq 2n-s$, where
$H^k_{\delta}({\frak g},\omega)$ denotes the cohomology of
$(\bigwedge^*_{\delta}({\frak g}^*, \omega)=\{ \alpha\in
\bigwedge^*({\frak g}^*) \mid \delta\alpha=0 \},\ d)$.
\end{corollary}

Notice that Theorems~\ref{gradosddelta} and~\ref{mapi} imply that
if two symplectic forms $\omega$ and $\omega'$ are cohomologous
then the $d\delta$-lemma holds up to degree $s$ for $\omega$ if
and only if it does for $\omega'$, and the map $i$ given in
(\ref{eqn:v1aa}) is bijective for all $k\geq 2n-s$ for $\omega$ if
and only if it is so for $\omega'$. Therefore, we can consider
symplectic forms up to cohomology class.

Next we consider an arbitrary symplectic form on some examples of
compact completely solvable manifolds, where we show explicit
calculations.

\begin{example} The Kodaira-Thurston manifold.
{\rm Let $G$ be the connected nilpotent Lie group of dimension $4$
given by $G=H\times \RR$, where $H$ is the Heisenberg group, that
is, the Lie group consisting of matrices of the form
 $$
 g=\left( \begin{array}{ccc} 1&x&z\\ 0&1&y\\ 0&0&1 \end{array} \right) ,
 $$
where $x,y,z \in {\RR}$. If $\Gamma'$ denotes the discrete
subgroup of $H$ consisting of matrices whose entries $x$, $y$ and
$z$ are integer numbers, then the quotient space
$KT=\Gamma{\backslash} G$, where $\Gamma=\Gamma'\times \ZZ$, is a
compact manifold.

A global system of coordinates $(x,y,z)$ for $H$ is given by
$x(g)=x$, $y(g)=y$, $z(g)=z$, and a standard calculation shows
that a basis for the left invariant $1$--forms on $H$ consists of
$\{dx, dy, dz-xdy\}$. Thus, if $t$ denotes the standard coordinate
for $\RR$, then $\{ \alpha=-dx, \beta=dy, \gamma=dt, \tau=dz-xdy
\}$ is a basis of the dual ${\frak g}^*$ of the Lie algebra
${\frak g}$ of $G$ with Chevalley-Eilenberg differential given by
 $$
 d\alpha=d\beta=d\gamma=0, \quad
 d\tau=\alpha \wedge \beta.
 $$
So, the Chevalley-Eilenberg cohomology of ${\frak g}$ is given by
\begin{eqnarray*}
 H^0({\frak g}) &=& \la 1\ra, \\
 H^1({\frak g}) &=& \la [\alpha], [\beta], [\gamma]\ra,\\
 H^2({\frak g}) &=& \la
   [\alpha \wedge \gamma], [\alpha \wedge \tau], [\beta \wedge \gamma], [\beta\wedge \tau]\ra,\\
 H^3({\frak g}) &=& \la [\alpha \wedge \beta \wedge \tau], [\alpha \wedge \gamma \wedge \tau],
   [\beta\wedge \gamma \wedge \tau]\ra,\\
 H^4({\frak g}) &=& \la [\alpha \wedge \beta \wedge \gamma \wedge \tau]\ra.
\end{eqnarray*}

For any element $\omega\in \bigwedge^2({\frak g}^*)$ satisfying
$d\omega=0$ there exists $a,b,c,e\in \RR$ such that
 $$
 [\omega]= a\, [\alpha \wedge \gamma] + b\, [\beta \wedge \gamma] +
 c\, [\alpha \wedge \tau] + e\, [\beta\wedge \tau].
 $$
Since $[\omega]^2= 2(bc-ae) [\alpha \wedge \beta\wedge\gamma\wedge
\tau]$, we conclude that $[\omega]^2\not=0$ if and only if $ae
\not= bc$. Hence, up to cohomology class, we can consider that any
symplectic form on ${\frak g}$ is given by
\begin{equation}\label{gen1-omega}
\omega= a\, \alpha \wedge \gamma + b\, \beta \wedge \gamma + c\,
\alpha \wedge \tau + e\, \beta\wedge \tau, \quad\quad ae-bc
\not=0.
\end{equation}
Moreover, notice that for the new basis of ${\frak g}^*$ given by
 $$
 \alpha'=(ae-bc)(a\, \alpha + b\, \beta),\quad
 \beta'=\frac{1}{ae-bc}(c\, \alpha + e\, \beta),\quad
 \gamma'=\frac{1}{ae-bc}\,\gamma,\quad \tau'=(ae-bc) \tau,
 $$
the differential $d$ expressed again as
 $$
 d\alpha'=d\beta'=d\gamma'=0,\quad
 d\tau'=\alpha' \wedge \beta'.
 $$
Now, with respect to this basis the symplectic
form~(\ref{gen1-omega}) is given by
 $$
 \omega= \alpha' \wedge \gamma' + \beta' \wedge \tau'.
 $$
Therefore, we can suppose without loss of generality that $a=e=1$
and $b=c=0$ in~(\ref{gen1-omega}).

Observe that $[\omega] \cup [\beta]=0$ in $H^3({\frak g})$, and
$\dim H^3_{\rm hr}({\frak g},\omega)=2<3=\dim H^3({\frak g})$. It
follows from Corollary~\ref{Nomizu-type-th} that for any
symplectic form $\omega$ on $KT$, the compact symplectic manifold
$(KT,\omega)$ is not $1$--Lefschetz, and $H^k_{\rm hr}(KT,\omega)
= H^k(KT)$ for $k\not=3$, but $\dim H^3_{\rm hr}(KT,\omega)=2 <
3=b_3(KT)$. Notice that any non-toral compact symplectic
nilmanifold $(M=\G\backslash G, \omega)$ is $0$--Lefschetz but not
$1$--Lefschetz~\cite{BG}.

We study next the $d\delta$-lemma for any symplectic form $\omega$
on the Kodaira-Thurston manifold. By
Corollary~\ref{d-delta-invariant}, the $d\delta$-lemma is
satisfied up to degree $s=0$ if and only if it is satisfied at the
level of the Lie algebra ${\frak g}$. Let us denote by $\{ X,Y,Z,T
\}$ the basis of ${\frak g}$ dual to $\{ \alpha, \beta, \gamma,
\tau \}$, and let $\omega$ be a symplectic form on ${\frak g}$
given by~(\ref{gen1-omega}) with $a=e=1$ and $b=c=0$. Then, the
isomorphism $\natural\colon {\frak g} \longrightarrow {\frak g}^*$
is given by
 $$
 \natural(X)=\gamma, \ \quad \natural(Y)=\tau, \ \quad
 \natural(Z)=-\alpha, \ \quad \natural(T)=-\beta.
 $$
Therefore,
 $$
 G =-\natural^{-1}(\omega) = - X\wedge Z - Y\wedge T.
 $$

In degree $1$, we must determine the spaces
$\delta(\bigwedge^2({\frak g}^*))\cap\ker d$,
$d(\bigwedge^0({\frak g}^*))\cap\ker \delta$ and
$d\delta(\bigwedge^1({\frak g}^*))$. Notice that
$\delta(\bigwedge^1({\frak g}^*))\subset \bigwedge^0({\frak g}^*)
= \RR$ and $d(\bigwedge^0({\frak g}^*))=\{0\}$. Using that
$\delta\mu=i_G(d\mu)$ for any $\mu\in \bigwedge^2({\frak g}^*)$,
an easy calculation shows that $\delta(\bigwedge^2({\frak
g}^*))=\la \beta \ra$, in fact $\beta=\delta(-\gamma\wedge\tau)$.
Since $\ker d=\langle \alpha,\beta,\gamma \rangle$, we have
 $$
 \delta(\bigwedge\phantom{i}\!\!\!^2({\frak g}^*))\cap \ker d=\la
 \beta \ra \not=\{0\} = d\delta(\bigwedge\phantom{i}\!\!\!^1({\frak g}^*)),
 $$
and the $d\delta$-lemma is not satisfied in degree $1$. Moreover,
$d(\bigwedge^1({\frak g}^*))\cap \ker \delta=\la \alpha\wedge\beta
\ra \not= \{0\} = d\delta(\bigwedge^2({\frak g}^*))$.

Applying the symplectic star operator, we get that the element
$\alpha\wedge\beta\wedge\gamma=*\beta \in d(\bigwedge^2({\frak
g}^*))\cap\ker \delta$, but it does not belong to the space
$d\delta(\bigwedge^3({\frak g}^*))=\{0\}$.

Therefore, for any symplectic form $\omega$ on $KT$ the
$d\delta$--lemma is satisfied only up to degree~$0$, according to
Theorem~\ref{gradosddelta}.

Notice that in general for any symplectic form on a nilpotent Lie
algebra the map $L^{n-1}$ is never injective~\cite{BG}. From
Theorem~\ref{gradosddelta} it follows that $\Im d\cap  \ker
\delta=\Im d\delta$ on $\bigwedge^1({\frak g}^*)$ and
$\Im\delta\cap\ker d =\Im d\delta$ on $\bigwedge^{2n-1}({\frak
g}^*)$, in fact these spaces are all zero, but either
$\Im\delta\cap \ker d =\Im d\delta$ fails on $\bigwedge^1({\frak
g}^*)$ or $\Im d\cap \ker\delta  =\Im d\delta$ fails on
$\bigwedge^2({\frak g}^*)$. By duality, $\Im d\cap \ker \delta=\Im
d\delta$ fails on $\bigwedge^{2n-1}({\frak g}^*)$, or
$\Im\delta\cap \ker d  =\Im d\delta$ fails on
$\bigwedge^{2n-2}({\frak g}^*)$.

Finally, we study the cohomology $H^*_{\delta}$. At the level of
${\frak g}$, the cohomology groups $H^k_{\delta}({\frak
g},\omega)$ are:
\begin{eqnarray*}
 H^0_{\delta}({\frak g},\omega) &=& \la 1\ra, \\
 H^1_{\delta}({\frak g},\omega) &=& \la [\alpha], [\beta], [\gamma]\ra,\\
 H^2_{\delta}({\frak g},\omega) &=& \la
   [\alpha \wedge \gamma], [\alpha \wedge \tau], [\beta \wedge \gamma], [\beta\wedge \tau]\ra,\\
 H^3_{\delta}({\frak g},\omega) &=& \la [\alpha \wedge \beta \wedge \gamma], [\alpha \wedge \gamma \wedge \tau],
   [\beta\wedge \gamma \wedge \tau]\ra,\\
 H^4_{\delta}({\frak g},\omega) &=& \la [\alpha \wedge \beta \wedge \gamma \wedge \tau]\ra.
\end{eqnarray*}
Therefore, $i\colon H^k_{\delta}({\frak g},\omega) \longrightarrow
H^k({\frak g})$ is bijective for all $k\not= 3$, because
$i([\alpha \wedge \beta \wedge \gamma])=0$ in $H^3({\frak g})$, in
fact $\alpha \wedge \beta \wedge \gamma=d(-\gamma\wedge\tau)$.
>From Corollary~\ref{mapi-invariant} we have that for any
symplectic form $\omega$ on $KT$ the map $i\colon
H^k_{\delta}(KT,\omega) \longrightarrow H^k(KT)$ is bijective for
$k=4$, but not for $k=3$, according to Theorem~\ref{mapi}. }
\end{example}

\begin{example} A six-dimensional solvmanifold.
{\rm  Let $G$ be the connected completely solvable Lie group of
dimension $6$ consisting of matrices of the form
 $$
 g=\pmatrix{e^t&0&xe^t&0&0&y_1 \cr
 0&e^{-t}&0&xe^{-t}&0&y_2\cr 0&0&e^t&0&0&z_1\cr 0&0&0&e^{-t}&0&z_2
 \cr 0&0&0&0&1&t \cr 0&0&0&0&0&1 \cr},
 $$
where $t, x, y_i, z_i \in \RR$ ($i=1,2$). The Lie group $G$ has a
discrete subgroup $\G$ such that the quotient space $M
=\Gamma\backslash G$ is compact~\cite{FLS}.

A global system of coordinates $(t,x,y_1,y_2,z_1,z_2)$ for $G$ is
defined by $t(g)=t$, $x(g)=x$, $y_i(g)=y_i$, $z_i(g)=z_i$, and a
standard calculation shows that a basis for the left invariant
$1$--forms on $G$ consists of
    $$
    \{\alpha=dt,\ \beta=dx,\  \gamma_1=e^{-t}dy_1-xe^{-t}dz_1,\
    \gamma_2=e^tdy_2-xe^tdz_2,\ \tau_1=e^{-t}dz_1,\ \tau_2=e^{t}dz_2\}.
    $$
Hence, $\{ \alpha, \beta, \gamma_1, \gamma_2, \tau_1, \tau_2 \}$
is a basis of the dual ${\frak g}^*$ of the Lie algebra ${\frak
g}$ of $G$ with Chevalley-Eilenberg differential given by
 $$
 d\alpha=d\beta=0, \quad
 d\gamma_1=-\alpha \wedge \gamma_1 - \beta \wedge \tau_1,\quad
 d\gamma_2=\alpha \wedge \gamma_2 - \beta \wedge \tau_2,\quad
 d\tau_1=- \alpha \wedge \tau_1,\quad
 d\tau_2= \alpha \wedge \tau_2.
 $$
Now, a direct calculation shows that the Chevalley-Eilenberg
cohomology of ${\frak g}$ is given by
 \begin{eqnarray*}
 H^0({\frak g}) &=& \la 1\ra,\\
 H^1({\frak g}) &=& \la [\alpha], [\beta]\ra,\\
 H^2({\frak g}) &=& \la [\alpha \wedge \beta],
 [\gamma_1 \wedge \tau_2 + \gamma_2 \wedge \tau_1], [\tau_1\wedge \tau_2]\ra,\\
 H^3({\frak g}) &=& \la
 [\alpha\wedge (\gamma_1 \wedge \tau_2 + \gamma_2 \wedge\tau_1)], [\alpha \wedge \tau_1 \wedge \tau_2],
  [\beta\wedge \gamma_1 \wedge \gamma_2],
  [\beta\wedge (\gamma_1 \wedge \tau_2 + \gamma_2 \wedge \tau_1)]\ra,\\
 H^4({\frak g}) &=& \la [\alpha \wedge \beta \wedge \gamma_1 \wedge \gamma_2],
   [\alpha \wedge \beta \wedge (\gamma_1 \wedge \tau_2 + \gamma_2 \wedge \tau_1)],
   [\gamma_1 \wedge \gamma_2 \wedge \tau_1\wedge \tau_2]\ra,\\
 H^5({\frak g}) &=& \la [\alpha \wedge \gamma_1 \wedge \gamma_2 \wedge \tau_1\wedge
   \tau_2], [\beta \wedge \gamma_1 \wedge \gamma_2 \wedge \tau_1\wedge
   \tau_2]\ra, \\
 H^6({\frak g}) &=& \la [\alpha \wedge \beta \wedge \gamma_1 \wedge \gamma_2 \wedge
       \tau_1\wedge \tau_2]\ra.
 \end{eqnarray*}

For any element $\omega\in \bigwedge^2({\frak g}^*)$ satisfying
$d\omega=0$ there exists $a,b,c\in \RR$ such that
 $$
 [\omega]= a\, [\alpha \wedge \beta] + b\,
 [\gamma_1 \wedge \tau_2 + \gamma_2 \wedge \tau_1] + c\, [\tau_1\wedge \tau_2].
 $$
Since $[\omega]^3=6ab^2[\alpha \wedge \beta\wedge\gamma_1 \wedge
\gamma_2 \wedge \tau_1 \wedge \tau_2]$, we conclude that
$[\omega]^3\not=0$ if and only if $ab\not=0$.

Thus, up to cohomology class, we can consider that any symplectic
form on ${\frak g}$ is given by
 \begin{equation}\label{gen2-omega}
 \omega= a\,\alpha \wedge \beta + b\, \gamma_1 \wedge \tau_2 + b\,
 \gamma_2 \wedge \tau_1 + c\, \tau_1\wedge \tau_2, \quad\quad
 a,b\not=0.
 \end{equation}
Let us consider the new basis of ${\frak g}^*$ given by
 $$
 \alpha'=\alpha,\,\, \beta'=a \, \beta,\,\, \gamma_1'=\sqrt{\frac{a}{b}}\,\left( b\gamma_1
 +\frac{c}{2}\,\tau_1\right),\,\,  \gamma_2'=\sqrt{\frac{a}{b}}\,\left( b\gamma_2
 -\frac{c}{2}\,\tau_2\right),\,\,  \tau'_1=\sqrt{\frac{b}{a}}\,\tau_1, \,\,
 \tau'_2=\sqrt{\frac{b}{a}}\,\tau_2,
 $$
if $ab>0$, or
 $$
 \alpha'=-\alpha, \,\, \beta'=-a\, \beta,\,\, \gamma_1'=\sqrt{-\frac{a}{b}}\, \left(b\gamma_2
 -\frac{c}{2}\,\tau_2\right),\qquad\qquad\qquad\qquad\qquad\qquad\qquad\qquad\qquad\qquad
 $$
 $$
 \qquad\qquad\qquad\qquad\qquad\qquad\qquad\qquad\gamma_2'=\sqrt{-\frac{a}{b}}\, \left(b\gamma_1
 +\frac{c}{2}\,\tau_1\right),\,\, \tau'_1= \sqrt{-\frac{b}{a}}\,\tau_2, \,\,
 \tau'_2=\sqrt{-\frac{b}{a}}\,\tau_1,
 $$
if $ab<0$. The differential $d$ also expressed as
 $$
 d\alpha'=d\beta'=0, \quad
 d\gamma'_1=-\alpha' \wedge \gamma'_1 - \beta' \wedge \tau'_1,\quad
 d\gamma'_2=\alpha' \wedge \gamma'_2 - \beta' \wedge \tau'_2,\quad
 d\tau'_1=- \alpha' \wedge \tau'_1,\quad
 d\tau'_2= \alpha' \wedge \tau'_2.
 $$
With respect to this basis the symplectic form~(\ref{gen2-omega})
is given by
 $$
 \omega= \alpha' \wedge \beta' + \gamma'_1
 \wedge \tau'_2 +  \gamma'_2 \wedge \tau'_1,
 $$
so we can suppose without loss of generality that $a=b=1$ and
$c=0$ in~(\ref{gen2-omega}).

Observe that $[\omega] \cup [\tau_1 \wedge \tau_2]=0$ in
$H^4({\frak g})$, but a simple computation shows that the product
by $[\omega]^2$ is an isomorphism between $H^1({\frak g})$ and
$H^5({\frak g})$. Moreover, $\dim H^4_{\rm hr}({\frak
g},\omega)=2<3=\dim H^4({\frak g})$. Therefore, for any symplectic
form $\omega$ on $M$, the compact symplectic manifold $(M,\omega)$
is $1$--Lefschetz, but not $2$--Lefschetz, and
Corollary~\ref{Nomizu-type-th} implies that $\dim H^k_{\rm
hr}(M,\omega) = b_k(M)$ for $k\not=4$, but $\dim H^4_{\rm
hr}(M,\omega)=2 < 3=b_4(M)$.

Next we study the $d\delta$-lemma for any symplectic form on the
compact solvmanifold $M$. Corollary~\ref{d-delta-invariant}
implies that the $d\delta$-lemma is satisfied up to degree $1$ on
$M$ if and only if it is satisfied on~${\frak g}$. Let us denote
by $\{ X,Y,Z_1,Z_2,T_1,T_2 \}$ the basis of ${\frak g}$ dual to
$\{ \alpha, \beta, \gamma_1, \gamma_2, \tau_1, \tau_2 \}$, and let
$\omega$ be a symplectic form on ${\frak g}$ given
by~(\ref{gen2-omega}) with $a=1$, $b=1$ and $c=0$. Then, the
isomorphism $\natural\colon {\frak g} \longrightarrow {\frak g}^*$
is given by
 $$
 \natural(X)=\beta, \ \quad \natural(Y)=-\alpha, \ \quad
 \natural(Z_1)=\tau_2, \ \quad \natural(Z_2)=\tau_1, \ \quad
 \natural(T_1)=-\gamma_2, \ \quad \natural(T_2)=-\gamma_1.
 $$
Therefore, $G=-\natural^{-1}(\omega)$ is given by
 $$
 G = - X\wedge Y -Z_1\wedge T_2 - Z_2\wedge T_1.
 $$

In degree $1$ we must consider the spaces
$\delta(\bigwedge^2({\frak g}^*))\cap\ker d$,
$d(\bigwedge^0({\frak g}^*))\cap\ker \delta$ and
$d\delta(\bigwedge^1({\frak g}^*))$. Since
$\delta(\bigwedge^1({\frak g}^*))\subset \bigwedge^0({\frak g}^*)
= \RR$, the $d\delta$-lemma is satisfied in degree $1$ if and only
if
 $$
 \delta(\bigwedge\phantom{i}\!\!\!^2({\frak g}^*))\cap\ker d =\{ 0
 \}.
 $$
Using that $\delta\mu=i_G(d\mu)$ for any $\mu\in
\bigwedge^2({\frak g}^*)$, a direct calculation shows that the
space $\delta(\bigwedge^2({\frak g}^*))$ is generated by
$\gamma_1$, $\gamma_2$, $\tau_1$ and $\tau_2$. Since $\ker
d=\langle \alpha,\beta\rangle$, the $d\delta$-lemma holds in
degree~$1$.

In degree $2$ we must compare the spaces
$\delta(\bigwedge^3({\frak g}^*))\cap\ker d$,
$d(\bigwedge^1({\frak g}^*))\cap\ker \delta$ and
$d\delta(\bigwedge^2({\frak g}^*))$. It is easy to check that
$d(\bigwedge^1({\frak g}^*))\subset \ker \{\delta\colon
\bigwedge^2({\frak g}^*)\longrightarrow \bigwedge^1({\frak
g}^*)\}$. Therefore,
 $$
 d(\bigwedge\phantom{i}\!\!\!^1({\frak g}^*))\cap\ker \delta
 =d(\bigwedge\phantom{i}\!\!\!^1({\frak g}^*))=\langle
 d\gamma_1,d\gamma_2,d\tau_1,d\tau_2 \rangle=d\delta(\bigwedge\phantom{i}\!\!\!^2({\frak
 g}^*)),
 $$
so this space is generated by $\alpha \wedge \gamma_1 +
\beta\wedge \tau_1$, $\alpha \wedge \gamma_2 - \beta \wedge
\tau_2$, $\alpha \wedge \tau_1$ and $\alpha \wedge \tau_2$.

However, a long but direct calculation shows that
 $$
 \delta(\bigwedge\phantom{i}\!\!\!^3({\frak g}^*))\cap\ker d= \langle
 d\gamma_1,d\gamma_2,d\tau_1,d\tau_2,\tau_1\wedge\tau_2 \rangle\
 \not\subset\ d(\bigwedge\phantom{i}\!\!\!^1({\frak
 g}^*))\cap \ker \delta.
 $$
In fact, notice that
 $$
 \begin{array}{rl}
 \delta(\alpha\wedge\gamma_1\wedge\tau_2)
 & = i_Gd(\alpha\wedge\gamma_1\wedge\tau_2) - \, di_G(\alpha\wedge\gamma_1\wedge\tau_2) \\[6pt]
 & = i_G(\alpha\wedge\beta\wedge\tau_1\wedge\tau_2) +
 d(\alpha)\\[6pt]
 & = - \tau_1\wedge\tau_2,
 \end{array}
 $$
and $d(\tau_1\wedge\tau_2)=0$. Thus, the $d\delta$-lemma is
satisfied up to degree~$1$, but it does not hold up to degree~$2$.
Therefore, for any symplectic form $\omega$ on $M$ the
$d\delta$--lemma is satisfied only up to degree~$1$, according to
Theorem~\ref{gradosddelta}.

Notice that the element
$\alpha\wedge\beta\wedge\tau_1\wedge\tau_2=*(-\tau_1\wedge\tau_2)
\in d(\bigwedge^3({\frak g}^*))\cap \ker \delta$ does not belong
to the space $d\delta(\bigwedge^4({\frak g}^*))$.

\medskip

Finally, the cohomology groups $H^k_{\delta}({\frak g},\omega)$
are given by:
\begin{eqnarray*}
 H^0_{\delta}({\frak g},\omega) &=& \la 1\ra,\\
 H^1_{\delta}({\frak g},\omega) &=& \la [\alpha], [\beta]\ra,\\
 H^2_{\delta}({\frak g},\omega) &=& \la [\alpha \wedge \beta],
 [\gamma_1 \wedge \tau_2 + \gamma_2 \wedge \tau_1], [\tau_1\wedge \tau_2]\ra,\\
 H^3_{\delta}({\frak g},\omega) &=& \la
 [\alpha\wedge (\gamma_1 \wedge \tau_2 + \gamma_2 \wedge
  \tau_1)], [\alpha \wedge \tau_1 \wedge \tau_2],
  [\beta\wedge \gamma_1 \wedge \gamma_2],
  [\beta\wedge (\gamma_1 \wedge \tau_2 + \gamma_2 \wedge \tau_1)]\ra,\\
 H^4_{\delta}({\frak g},\omega) &=& \la
   [\alpha \wedge \beta \wedge (\gamma_1 \wedge \tau_2 + \gamma_2 \wedge
   \tau_1)],
   [\alpha \wedge \beta \wedge \tau_1 \wedge \tau_2],
   [\gamma_1 \wedge \gamma_2 \wedge \tau_1\wedge \tau_2]\ra,\\
 H^5_{\delta}({\frak g},\omega) &=& \la [\alpha \wedge \gamma_1 \wedge \gamma_2 \wedge \tau_1\wedge
   \tau_2], [\beta \wedge \gamma_1 \wedge \gamma_2 \wedge \tau_1\wedge
   \tau_2]\ra, \\
 H^6_{\delta}({\frak g},\omega) &=& \la [\alpha \wedge \beta \wedge \gamma_1 \wedge \gamma_2 \wedge
       \tau_1\wedge \tau_2]\ra.
\end{eqnarray*}
Thus, the map $i\colon H^k_{\delta}({\frak g},\omega)
\longrightarrow H^k({\frak g})$ is bijective for all $k\not= 4$.
In fact, since $d(\alpha\wedge\gamma_1\wedge\tau_2)=\alpha \wedge
\beta \wedge \tau_1 \wedge \tau_2$ we have that $i([\alpha \wedge
\beta \wedge \tau_1 \wedge \tau_2])=0$ in $H^4({\frak g})$. By
Corollary~\ref{mapi-invariant} we conclude that for any symplectic
form $\omega$ on $M$ the map $i\colon H^k_{\delta}(M,\omega)
\longrightarrow H^k(M)$ is bijective for $k=5,6$, but not for
$k=4$, according to Theorem~\ref{mapi}.

}
\end{example}

\bigskip

\noindent {\bf Acknowledgments.} We thank to G. Cavalcanti for
useful conversations. This work has been partially supported
through grants MCyT (Spain) Project BFM2001-3778-C03-02/03, UPV
00127.310-E-14813/2002 and MTM2004-07090-C03-01.

{\small

\vspace{0.15cm}

\noindent{\sf M. Fern\'andez:} Departamento de Matem\'aticas,
Facultad de Ciencia y Tecnolog\'{\i}a, Universidad del Pa\'{\i}s
Vasco, Apartado 644, 48080 Bilbao, Spain. {\sl E-mail:}
mtpferol@lg.ehu.es

\vspace{0.15cm}

\noindent{\sf V. Mu\~noz:} Instituto de Matem\'aticas y
F\'{\i}sica Fundamental, Consejo Superior de Investigaciones
Cient\'{\i}ficas, Serrano 113 bis, 28006 Madrid, Spain. {\sl
E-mail:} vicente.munoz@imaff.cfmac.csic.es

\vspace{0.15cm}

\noindent{\sf L. Ugarte:} Departamento de Matem\'aticas, Facultad
de Ciencias, Universidad de Zaragoza, Campus Plaza San Francisco,
50009 Zaragoza, Spain. {\sl E-mail:} ugarte@unizar.es}

\label{endxyzt}

\begin{thebibliography}{33}

\bibitem{BG}
C. Benson and C.S. Gordon, K\"ahler and symplectic structures on
nilmanifolds, {\em Topology\/} {\bf 27} (1988), 513--518.

 \bibitem{Bry}
J.L. Brylinski, A differential complex for Poisson manifolds,
     {\em J. Diff. Geom.\/}   {\bf 28} (1988),  93--114.

\bibitem{Cav1}
G.R. Cavalcanti, The Lefschetz property, formality and blowing up
in symplectic geometry, Preprint {\tt math.SG/0403067}.

\bibitem{Cav2}
G.R. Cavalcanti, New  aspects of the $dd^c$-lemma, Ph. D. Thesis,
University of Oxford, $2004$.

\bibitem{FLS}
 M. Fern\'andez, M. de Le\'on and M. Saralegui, A six dimensional
compact symplectic solvmanifold without K\"ahler structures,
 {\em Osaka J. Math.\/} {\bf 33} (1996), 19--35.

\bibitem{FM}
M. Fern\'andez and V. Mu\~noz, Formality of Donaldson
submanifolds,
 {\em Math.\ Zeit.,\/} To appear.

\bibitem{FMU}
M. Fern\'andez, V. Mu\~noz and L. Ugarte, Weakly Lefschetz
symplectic manifolds, Preprint {\tt math.SG/0404479}.

\bibitem{Gui}
V. Guillemin, Symplectic Hodge theory and the $d\delta$--lemma,
Preprint, Massachusets Institute of Technology, 2001.

\bibitem{Hat}
A. Hattori, Spectral sequence in the de Rham cohomology of fibre
bundles,
 {\em J. Fac. Sci.  Univ. Tokyo\/} {\bf 8} (1960), 298--331.

\bibitem{IRTU}
R. Ib\'{a}\~{n}ez, Y. Rudyak, A. Tralle and L. Ugarte, On
symplectically harmonic forms on $6$--dimensional nilmanifolds,
{\em Comment. Math. Helv.\/} {\bf 76} (2001), 89--109.


\bibitem{Kos}
J.L. Koszul, Crochet de Schouten-Nijenhuis et cohomologie, in {\em
Elie Cartan et les  Math. d'Aujour d'Hui}, Ast\'{e}risque
hors-s\'{e}rie (1985), 251--271.

\bibitem{Lib}
P. Libermann, Sur le probl\`eme d'equivalence de certaines
structures infinitesimales regulieres, {\em Ann. Mat. Pura
Appl.\/} {\bf 36} (1954), 27--120.

\bibitem{LM}
P. Libermann and C. Marle, {\em Symplectic Geometry and Analytical
Mechanics},  Kluwer, Dordrecht, 1987.

\bibitem{Lich}
A. Lichnerowicz, Les vari\'{e}t\'{e}s de Poisson et les
alg\'{e}bres de Lie associ\'{e}es, { \em  J. Diff. Geom.\/} {\bf
12} (1977), 253--300.

\bibitem{LS}
Y. Lin and R. Sjamaar, Equivariant symplectic Hodge theory and the
$d_{G}\delta$--lemma, {\em J. Symplectic Geom.\/}, To appear,
Preprint {\tt math.SG/0310048}.


\bibitem{Mat}
O. Mathieu, Harmonic cohomology classes of symplectic manifolds,
  {\em Comment. Math. Helv.\/} {\bf 70} (1995), 1--9.


\bibitem{Mer}
S. Merkulov, Formality of canonical symplectic complexes and
Frobenius manifolds,
   {\em Internat. Math. Res. Notices\/} {\bf 14} (1998), 723--733.

\bibitem{No}
K. Nomizu, On the cohomology of compact homogeneous spaces of
nilpotent Lie groups, {\em Ann. of Math.\/} {\bf 59} (1954),
531-538.

\bibitem{Tu}
W.P. Thurston, Some simple examples of symplectic manifolds,
  {\em Proc. Amer.  Math. Soc.\/} {\bf 55} (1976), 467--468.

\bibitem{wells}
 R. Wells, Differential analysis on complex manifolds. Second edition.
 Graduate Texts in Mathematics {\bf 65}. Springer-Verlag, New
 York-Berlin, 1980.

\bibitem{Yam}
 T. Yamada, Harmonic cohomology groups of compact symplectic nilmanifolds,
  {\em Osaka J. Math.\/} {\bf 39} (2002), 363--381.

\bibitem{Ya}
 D. Yan, Hodge structure on symplectic manifolds,
  {\em Adv. Math.\/} {\bf 120} (1996), 143--154.


\end{thebibliography}
\end{document}